\documentclass[12pt,reqno]{amsart}
\usepackage{amssymb,amsmath,amscd}
\usepackage[english]{babel}
\makeatother
\usepackage{wrapfig}

\newcommand{\lit}[3]{\vspace*{0.7mm}\par\noindent\makebox[5.2mm][r]{#1.~}\parbox[t]{159.8mm}{{\textit{#2}}~{#3}}\hspace*{-1.6mm}}
\author{A. A. Kon'kov}
\address{Department of Differential Equations,
Faculty of Mechanics and Mathematics,
Mo\-s\-cow Lo\-mo\-no\-sov State University,
Vorobyovy Gory,
119992 Moscow, Russia}
\email{konkov@mech.math.msu.su}
\title[On the absence of solutions]{On the absence of solutions of differential inequalities with $\infty$-Laplacian}
\date{}

\begin{document}



%
%
%




\begin{abstract}\noindent
For differential inequalities with the $\infty$-Laplacian in the principal part, we obtain conditions for the absence of solutions in unbounded domains. Examples are given to demonstrate the accuracy of these conditions.
\end{abstract}

\maketitle

\textbf{Introduction.} 
We study solutions of the problem
\begin{equation}
	\Delta_\infty
	u
	\ge
	f (x, u)
	\quad
	\mbox{in } \Omega,
	\qquad
	\left.
		u
	\right|_{
		\partial \Omega
	}
	=
	0,
	\label{1.1}
\end{equation}
where $\Omega$ is unbounded domain ${\mathbb R}^n$, $n \ge 2$,
$$
	\Delta_\infty
	u
	=
	\frac{1}{2}
	\nabla
	|\nabla u|^2
	\nabla u
	=
	\sum_{i,j=1}^n
	\frac{
		\partial^2 u
	}{
		\partial x_i \partial x_j
	}
	\frac{
		\partial u
	}{
		\partial x_i
	}
	\frac{
		\partial u
	}{
		\partial x_j
	}
$$
is the $\infty$-Laplacian, and  $f$ is some function such that
\begin{equation}
	f (x, t) > 0
	\label{1.2}
\end{equation}
for all $t > 0$ and $x \in \Omega$.
Let us agree to denote by $B_r$ and $S_r$, respectively, the open ball and sphere in ${\mathbb R}^n$ of radius $r > 0$ centered at zero. We shall assume that $S_r \cap \Omega \ne \emptyset$ for all $r > r_*$, where $r_* > 0$ is some real number.
Also let there be real numbers $\sigma > 1$ and $\theta > 1$ such that\begin{equation}
	\inf_{
		x \in \Omega_{r, \sigma},
		\,
		\zeta \in (t / \theta, t \theta)
	}
	f (x, |x| \zeta)
	\ge
	q (r)
	g (t)
	\label{1.5}
\end{equation}
for all $t > 0$ and $r > r_*$, where
$
	\Omega_{r, \sigma} 
	= 
	\{ 
		x \in \Omega 
		: 
		r / \sigma < |x| < r \sigma
	\}
$
and $q$ and $g$ are non-negative measurable functions such that
$$
	\inf_K g > 0
$$
for any compact set $K \subset (0, \infty)$.

Solutions of the problem~\eqref{1.1} is understood in the viscosity sense~[1, 2].
If $\partial \Omega = \emptyset$, i.e. $\Omega = {\mathbb R}^n$, then the boundary condition in~\eqref{1.1} is assumed to be satisfied automatically.

The phenomenon of the absence of non-trivial solutions for differential equations and inequalities (blow-up) has been studied by many authors~[3--19].
In the case when the main part of the differential operator has a divergent form, the method of nonlinear capacity has proven itself well. The main ideas of this method are presented in~[6, 18].
For inequalities~\eqref{1.1} this method obviously does not work.
Under the assumption that $\Omega$ is a bounded domain, the blow-up phenomenon for inequalities of the form~\eqref{1.1} was studied in~[20--22].
These papers dealt with the so-called ``boundary blow-up'' or, in other words, ``large solutions'', tending to infinity when approaching $\partial \Omega$.
The dependence of the right-hand side on the spatial variable was essentially not taken into account in the works~[20--22].
In the case of unbounded domains, the conditions for the absence of nontrivial solutions to the problem~\eqref{1.1} have not been known until now. Theorem~1 proved in our paper eliminates this gap.

{\bf Definition 1.}
A function 
$u : \overline{\Omega} \to {\mathbb R}$ 
is called upper semicontinuous at a point
$x \in \overline{\Omega}$
if
$$
	\limsup_{y \to x} u (y) \le u (x).
$$
We say that $u$ is upper semicontinuous on the set ${\mathcal G} \subset \overline{\Omega}$ if it is upper semicontinuous at every point of $x \in {\mathcal G}$.

{\bf Definition 2.}
The second-order superjet $J_\Omega^{2,+} u (x)$ of a function $u : \Omega \to {\mathbb R}$ at a point $x \in \Omega$ is the set of ordered pairs $(v, A )$, where $v \in {\mathbb R}^n$ is some vector and $A$ is a symmetric $n \times n$ matrix such that$$
	u (y) - u (x) 
	\le
	\langle v, y - x \rangle
	+
	\frac{1}{2}
	\langle A (y - x), y - x \rangle
	+
	\bar{\bar{o}} (|y - x|^2)
$$
for all $y$ from some neighborhood of the point $x$. Here and below, angle brackets denote the standard scalar product in ${\mathbb R}^n$.

Note that in the case of $u \in C^2 (\Omega)$ the set $J_\Omega^{2,+} u (x)$ is not empty for any point $x \in \Omega$, because
$$
	(D u (x), D^2 u (x)) 
	\in 
	J_\Omega^{2,+} u (x),
$$
where $D u (x)$ is the gradient vector and $D^2 u (x)$ is the matrix of second derivatives of the function $u$ at the point~$x$.

{\bf Definition 3.}
A function $u \ge 0$ is called a non-negative solution of problem~\eqref{1.1} if $u$ is upper semicontinuous on $\overline{\Omega}$, vanishes on $\partial \Omega$
and for any $x \in \Omega$ such that $J_\Omega^{2,+} u (x) \ne \emptyset$ the relation
\begin{equation}
	\langle A v, v \rangle
	\ge
	f (x, u (x))
	\label{d1.3.1}
\end{equation}
is valid for all $(v, A) \in J_\Omega^{2,+} u (x)$.

In the general case~[1], an upper semicontinuous function $u$ is called a viscosity solution of the inequality
$$
	F (x, u, D u, D^2 u) \le 0
	\quad
	\mbox{in } \Omega,
$$
where $F$ is some function if for any point $x \in \Omega$ such
that $J_\Omega^{2,+} u (x) \ne \emptyset$ the relation
\begin{equation}
	F (x, u (x), v, A) \le 0
	\label{1.3}
\end{equation}
is valid fro all $(v, A) \in J_\Omega^{2,+} u (x)$.
In particular, taking
\begin{equation}
	F (x, u, D u, D^2 u) 
	= 
	f (x, u) 
	- 
	\sum_{i,j=1}^n
	\frac{
		\partial^2 u
	}{
		\partial x_i \partial x_j
	}
	\frac{
		\partial u
	}{
		\partial x_i
	}
	\frac{
		\partial u
	}{
		\partial x_j
	}
	\label{1.4}
\end{equation}
in formula~\eqref{1.3}, we obtain~\eqref{d1.3.1}.

In order for the classical solutions to be simultaneously viscous, the function $F$ is subjected to the condition of degenerate ellipticity
$$
F (x, \zeta, v, A) \le F (x, \zeta, v, B)
$$
for all $x \in \Omega$, real numbers $\zeta$, vectors $v \in {\mathbb R}^n$, and symmetric ${n \times n}$-matrices $A \ge B$. The last inequality is understood in the sense of quadratic forms.
Namely, we say that $A \ge B$ if
$$
\langle A v, v \rangle
\ge
\langle B v, v \rangle
$$
for all $v \in {\mathbb R}^n$.
It is obvious that the function $F$ defined by~\eqref{1.4} satisfies  the condition of degenerate ellipticity.

\pagebreak

\textbf{1. Main results.} The following statements are valid.

\textbf{Theorem 1.} {\it 
Suppose that
\begin{equation}
	\int\limits_1^\infty
	(g (t) t)^{- 1 / 4}
	dt
	<
	\infty
	\quad
	\mbox{and}
	\quad
	\int\limits_{r_*}^\infty
	q (r)
	dr
	=
	\infty.
	\label{t2.1.1}
\end{equation}
Then any non-negative solution of~(\ref{1.1}) is identically equal to zero.
}

A special case of~\eqref{1.1} is the problem
\begin{equation}
	\Delta_\infty
	u
	\ge
	c (|x|) u^\lambda
	\quad
	\mbox{in } \Omega, 
	\qquad
	\left.
		u
	\right|_{
		\partial \Omega
	}
	=
	0,
	\label{2.1}
\end{equation}
where $c$ is a positive non-increasing function.

\textbf{Corollary 1.} {\it 
Suppose that $\lambda > 3$ and
$$
	\int\limits_{r_*}^\infty
	r^\lambda
	c (r)
	dr
	=
	\infty.
$$
Then any non-negative solution of~(\ref{2.1}) is identically equal to zero.
}

Theorem~1 and Corollary~1 will be proved in Paragraph~2.
Now we give examples de\-mon\-stra\-ting the accuracy of these results.

{\bf Example 1.} 
Consider the inequality
\begin{equation}
	\Delta_\infty u \ge c_0 (1 + |x|)^s u^\lambda
	\quad
	\mbox{in }
	{\mathbb R}^n,
	\label{e2.1.1}
\end{equation}
where $c_0 = const > 0$.
According to Corollary~1, the conditions
$$
	\lambda > 3
	\quad
	\mbox{and}
	\quad
	s \ge - \lambda - 1
$$
guarantee that any non-negative solution of~\eqref{e2.1.1} is identically equal to zero.
These conditions are exact. Indeed, it is not difficult to verify that, for
$\lambda > 3$ and $s < - \lambda - 1$, the expression
$$
	u (x)
	=
	\left\{
		\begin{aligned}
			&
			|x|^{- (4 + s) / (\lambda - 3)} - r_0^{- (4 + s) / (\lambda - 3)},
			&
			&
			|x| > r_0,
			\\
			&
			0,
			&
			&
			|x| \le r_0,
		\end{aligned}
	\right.
$$
where $r_0 > 0$ is some real number,
defines a non-negative solution of~\eqref{e2.1.1}.
At the same time, in the case of $\lambda \le 3$, by somewhat more complicated reasoning, one can show that the inequality~\eqref{e2.1.1} has a non-trivial non-negative solution for all $s \in {\mathbb R}$.

{\bf Example 2.} 
We examine the critical exponent $s = - \lambda - 1$ in the right-hand side of~\eqref{e2.1.1}.
Let $u \ge 0$ satisfy the inequality
\begin{equation}
	\Delta_\infty u 
	\ge 
	c_0 
	(1 + |x|)^{- \lambda - 1} 
	\ln^\mu (2 + |x|)
	u^\lambda
	\quad
	\mbox{in }
	{\mathbb R}^n,
	\label{e2.2.1}
\end{equation}
where $c_0 = const > 0$. Applying Corollary~1, we obtain that, for
\begin{equation}
	\lambda > 3
	\quad
	\mbox{and}
	\quad
	\mu \ge - 1,
	\label{e2.2.2}
\end{equation}
the function $u$ is identically equal to zero.
At the same time, by direct calculations, one can verify that, for $\lambda > 3$ and $\mu < -1$, the relation
$$
	u (x)
	=
	\left\{
		\begin{aligned}
			&
			|x| \ln^{- (\mu + 1) / (\lambda - 3)} |x| - r_0 \ln^{- (\mu + 1) / (\lambda - 3)} r_0,
			&
			&
			|x| > r_0,
			\\
			&
			0,
			&
			&
			|x| \le r_0,
		\end{aligned}
	\right.
$$
where $r_0 > 0$ is large enough, defines a non-negative solution~\eqref{e2.2.1}.
Thus, the second condition in~\eqref{e2.2.2} is exact. The first condition in~\eqref{e2.2.2} is also exact. It can be shown that, in the case of $\lambda \le 3$, inequality~\eqref{e2.2.1} has a non-trivial non-negative solution for all $\mu \in {\mathbb R}$.

{\bf Example 3.} 
Now we examine investigate the critical exponent $\lambda = 3$ in~\eqref{e2.1.1}.
Namely, consider the inequality
\begin{equation}
	\Delta_\infty u
	\ge
	c_0 (1 + |x|)^s
	u^3
	\ln^\nu (2 + u)
	\quad
	\mbox{in }
	{\mathbb R}^n,
	\label{e2.3.1}
\end{equation}
where $c_0 = const > 0$.
According to Theorem~1, if
\begin{equation}
	\nu > 4
	\quad
	\mbox{and}
	\quad
	s \ge - 4,
	\label{e2.3.2}
\end{equation}
then any non-negative solution of~\eqref{e2.3.1} is identically equal to zero.
In so doing, if $\nu > 4$ and $s < - 4$, then
$$
	u (x)
	=
	\left\{
		\begin{aligned}
			&
			e^{
				|x|^{
					(4 + s) / (4 - \nu)
				}
			}
			- 
			e^{
				r_0^{
					(4 + s) / (4 - \nu)
				}
			},
			&
			&
			|x| > r_0,
			\\
			&
			0,
			&
			&
			|x| \le r_0,
		\end{aligned}
	\right.
$$
where $r_0 > 0$ is large enough, is a non-negative solution of~\eqref{e2.3.1}.
If the first condition in~\eqref{e2.3.2} is violated, i.e. $\nu \le 4$, then~\eqref{e2.3.1} has a non-trivial non-negative solution for any $s \in {\mathbb R}$.
Thus, both conditions in~\eqref{e2.3.2} are exact.

\smallskip

\textbf{2. Proof of Theorem~1 and Corollary~1.} 
Below we assume that conditions~\eqref{t2.1.1} are satisfied.
Take a function $\omega \in C^\infty ({\mathbb R})$ positive on the interval $(-1, 1)$, equal to zero on the set ${\mathbb R} \setminus (-1, 1) $, and such that
$$
	\int\limits_{-\infty}^\infty
	\omega (t)
	dt
	=
	1.
$$
We denote
$$
	\omega_\delta (t) 
	= 
	\frac{1}{\delta}
	\omega 
	\left(
		\frac{t}{\delta}
	\right),
	\quad
	\delta > 0.
$$ 

\textbf{Lemma~1.} {\it 
Let $\mu > 1$, $\nu > 1$, and $0 < \alpha \le 1$ be some real numbers and, moreover,
$\eta : (0,\infty) \to (0,\infty)$ 
and
$H : (0,\infty) \to (0,\infty)$ be measuable functions satisfying the inequality
$$
	\eta (t)
	\le
	\inf_{
		(t / \mu, t \mu)
	}
	H
$$
for all $t \in (0, \infty)$.
Then
$$
	\left(
		\int\limits_{M_1}^{M_2}
		\eta^{-\alpha} (t)
		t^{\alpha - 1}
		\,
		dt
	\right)^{1 / \alpha}
	\ge
	C
	\int\limits_{M_1}^{M_2}
	\frac{
		d t
	}{
		H (t)
	}
$$
for all real numbers $M_1 > 0$ and $M_2 > 0$ such that $M_2 \ge \nu M_1$,
where the constant $C > 0$ depends only on $\alpha$, $\nu$ and $\mu$.
}

\textbf{Lemma~2.} {\it 
Let $0 < r_1 < r_2$ and $0 < \alpha \le 1$ be some real numbers and
$\eta : (r_1, r_2) \to [0, \infty)$ be measurable function. Then
$$
	\int\limits_{
		r_1
	}^{
		r_2
	}
	\left(
		\int\limits_{
			r_1
		}^\rho
		\eta (\xi)
		d\xi
	\right)^\alpha
	d\rho
	\ge
	C
	\int\limits_{
		r_1
	}^{
		r_2
	}
	\left(
		1
		-
		\frac{\xi}{r_2}
	\right)^\alpha
	(
		\xi \eta (\xi)
	)^\alpha
	d\xi,
$$
where the constant $C > 0$ depends only on $\alpha$.
}

\textbf{Lemma~3.} {\it 
Let
$0 < r_1 < r_2$
and
$0 < \alpha < 1$
be some real numbers and
$
	\eta : (r_1, r_2) \to [0, \infty)
$
be a measurable function. Then
\begin{equation}
	\left(
		\int\limits_{
			r_1
		}^{
			r_2
		}
		\eta (\xi)
		d\xi
	\right)^\alpha
	\ge
	C
	\int\limits_{
		r_1
	}^{
		r_2
	}
	\eta (\xi)
	\varkappa^{\alpha - 1} (\xi)
	d\xi,
	\label{l3.3.1}
\end{equation}
where
$$
	\varkappa (\xi)
	=
	\int\limits_\xi^{r_2}
	\eta (\zeta)
	d\zeta,
$$
and the constant $C > 0$ depends only on $\alpha$.
}

\textbf{Remark~1.}
If $\varkappa (\xi) = 0$ for some $\xi \in (r_1, r_2)$, then
$\eta (\zeta) = 0$ for almost all $\zeta \in (\xi, r_2)$.
In this case, we assume that
$
	\eta (\xi)
	\varkappa^{
		\alpha - 1
	}
	(\xi)
$
is also equal to zero on the right in~\eqref{l3.3.1}.

\textbf{Lemma~4.} {\it 
Let $\gamma > 0$, $\lambda > 1$, and $r_0 > r_*$ be some real numbers and, moreover,
$p : [r_*, \infty) \to [0, \infty)$ be a locally integrable function such that
$$
	\frac{
		r p (r)
	}{
		\int\limits_{r_*}^r
		p (\xi)
		d\xi
	}
	\le
	\gamma
$$
for all $r \ge r_0$.
Then
$$
	\int\limits_{r_*}^{
		\lambda r
	}
	p (\xi)
	d\xi
	\le
	C
	\int\limits_{r_*}^r
	p (\xi)
	d\xi
$$
for all $r \ge r_0$, where the constant $C > 0$ depends only on $\gamma$ and $\lambda$.
}

The proof of Lemma~1 can be found in~[13, Lemma~2.3].
Lemmas~2 and~3 were proved in~[12, Lemmas~2.2 and~2.3]. Lemma~4 is proved in~[14, Lemma~2.7].

\textbf{Lemma~5.} {\it 
Let $\mu > 1$ be some real number and $H : (0, \infty) \to (0, \infty)$ be a measurable function such that
$$
	\eta (t)
	=
	\inf_{
		(t / \mu, t \mu)
	}
	H
	>
	0
$$
for all $t \in (0, \infty)$ and, moreover,
\begin{equation}
	\int\limits_1^\infty
	(\eta (t) t)^{- 1 / 4}
	dt
	<
	\infty.
	\label{l3.5.2}
\end{equation}
Then there exists an infinitely smooth function $h : (0, \infty) \to (0, \infty)$ satisfying the following conditions:
\begin{enumerate}
\item[$1)$] 
$h (t) \le H (t)$ for all $t \in (0, \infty)$;
\item[$2)$] 
$h (t_2) \ge h (t_1)$ for all $t_2 \ge t_1 > 0$;
\item[$3)$] 
for any real number $\alpha > 0$ there is a constant $\beta > 0$,
depending only on $\alpha$ and $\mu$ such that
$h (\alpha t) \le \beta h (t)$ for all $t \ge 2$;
\item[$4)$] 
\begin{equation}
	\int\limits_1^\infty
	(h (t) t)^{- 1 / 4}
	dt
	<
	\infty.
	\label{l3.5.3}
\end{equation}
\end{enumerate}
}

{\bf Proof.} 
Take $t_k = \mu^{k/4}$, $k = 0,1,2,\ldots$.
Arguing by induction, we construct a sequence of monotonically nondecreasing functions
$h_k : [1, t_k] \to (0, \infty)$,
$k = 1,2,\ldots$,
with the following properties:
\begin{enumerate}
\item[i)]
$h_k (t) \le h_s (t) \le H (t)$
for all $t \in [1, t_s]$, $1 \le s \le k$;
\item[ii)]
$
	h_k (\mu^{1/4} t)
	\le
	\mu^2
	h_k (t)
$
for all $t \in [1, t_{k-1}]$, $k = 1,2,\ldots$;
\item[iii)]
\begin{equation}
	\int\limits_1^{t_k}
	(h_k (t) t)^{- 1 / 4}
	dt
	\le
	\frac{
		\mu^{1/2}
	}{
		1 - \mu^{-1/16}
	}
	\int\limits_1^{t_k}
	(\eta (t) t)^{- 1 / 4}
	dt
	\label{pl3.5.1}
\end{equation}
for all $k = 1,2,\ldots$.
\end{enumerate}
Denote
$$
	\gamma_k
	=
	\sup_{
		(t_{k-1}, t_{k+1})
	}
	\eta,
	\quad
	k = 1,2,\ldots.
$$
Let us put $h_1 = \gamma_1$.
Assume now that for some $k \ge 2$ the function $h_{k-1}$ is already known. We take
$$
	h_k (t)
	=
	\left\{
		\begin{aligned}
			&
			\min
			\{
				\gamma_k,
				\,
				h_{k-1} (t)
			\},
			&
			&
			t \in [1, t_{k-1}],
			\\
			&
			\min
			\{
				\gamma_k,
				\,
				h_{k-1} (t_{k-1})
			\},
			&
			&
			t \in [t_{k-1}, t_k]
		\end{aligned}
	\right.
$$
if
\begin{equation}
	\mu
	h_{k-1} (t_{k-1})
	\ge
	\gamma_k
	\label{pl3.5.2}
\end{equation}
and
$$
	h_k (t)
	=
	\left\{
		\begin{aligned}
			&
			h_{k-1} (t),
			&
			&
			t \in [1, t_{k-1}],
			\\
			&
			\frac{
				t_k - t
				+
				\mu
				(t - t_{k-1})
			}{
				t_k - t_{k-1}
			}
			h_{k-1} (t_{k-1}),
			&
			&
			t \in [t_{k-1}, t_k],
		\end{aligned}
	\right.
$$
if condition~\eqref{pl3.5.2} is not fulfilled.
It is easy to see that $h_k$ satisfies conditions i) and ii).
Let us show that $h_k$ satisfies inequality~\eqref{pl3.5.1}.
In the case $k = 1$ this inequality follows from the definitions of $h_1$ and $\gamma_1$.
Let $k \ge 2$ and 
\begin{equation}
	\int\limits_1^{t_s}
	(h_s (t) t)^{- 1 / 4}
	dt
	\le
	\frac{
		\mu^{1/2}
	}{
		1 - \mu^{-1/16}
	}
	\int\limits_1^{t_s}
	(\eta (t) t)^{- 1 / 4}
	dt
	\label{pl3.5.3}
\end{equation}
for all $1 \le s \le k - 1$.
At first, we assume that~\eqref{pl3.5.2} holds. Let $s \ge 0$ be the minimum integer such that $\mu h_k (t_s) \ge h_k (t_k)$.
It follows from the definition of the function $h_k$ that $h_s = h_k$ on the interval $[1, t_s]$.
Since
$$
	h_k (t) \ge \frac{h_k (t_k)}{\mu} \ge \frac{\gamma_k}{\mu^2}
$$ 
for all $t \in (t_s, t_k)$, we obtain
$$
	\int\limits_{t_s}^{t_k}
	(h_k (t) t)^{- 1 / 4}
	dt
	\le
	\frac{4}{3}
	\mu^{1 / 2}
	\gamma_k^{- 1 / 4}
	t_k^{3 / 4}.
$$
Combining the last inequality with the obvious relation
$$
	\int\limits_{
		t_{k-1}
	}^{
		t_k
	}
	(\eta (t) t)^{- 1 / 4}
	dt
	\ge
	\frac{4}{3}
	\gamma_k^{- 1 / 4}
	(t_k^{3 / 4} - t_{k - 1}^{3 / 4})
	=
	\frac{4}{3}
	\gamma_k^{- 1 / 4}
	(1 - \mu^{- 3 / 16})
	t_k^{3 / 4},
$$
we have
\begin{equation}
	\int\limits_{t_s}^{t_k}
	(h_k (t) t)^{- 1 / 4}
	dt
	\le
	\frac{
		\mu^{1/2}
	}{
		1 - \mu^{-3/16}
	}
	\int\limits_{
		t_{k-1}
	}^{
		t_k
	}
	(\eta (t) t)^{- 1 / 4}
	dt.
	\label{pl3.5.4}
\end{equation}
If $s = 0$, then~\eqref{pl3.5.4} immediately implies~\eqref{pl3.5.1}. Assume that $s \ge 1$.
In view of the fact that $s \le k - 1$, estimate~\eqref{pl3.5.3} is valid, summing which with~\eqref{pl3.5.4}, we derive~\eqref{pl3.5.1}.

Now assume that~\eqref{pl3.5.2} does not hold.
In this case, $h_k (t_k) = \mu h_k (t_{k-1})$.
Let $s \ge 1$ be the minimum integer such that
$h_k (t_k) = \mu^{k-s} h_k (t_s)$.
We have
$$
	\int\limits_{t_i}^{t_{i+1}}
	(h_k (t) t)^{- 1 / 4}
	dt
	=
	\mu^{ - (i - s) / 16}
	\int\limits_{t_s}^{t_{s+1}}
	(h_k (t) t)^{- 1 / 4}
	dt
$$
for all $s \le i \le k - 1$; therefore,
\begin{equation}
	\int\limits_{t_s}^{t_k}
	(h_k (t) t)^{- 1 / 4}
	dt
	\le
	\frac{
		1
	}{
		1 - \mu^{- 1 / 16}
	}
	\int\limits_{t_s}^{t_{s+1}}
	(h_k (t) t)^{- 1 / 4}
	dt.
	\label{pl3.5.5}
\end{equation}
At the same time, taking into account the monotonicity of the function $h_k$ and the fact that
$$
	h_k (t_s) 
	= 
	h_k (t_{s-1}) 
	\ge 
	\frac{\gamma_s}{\mu},
$$ 
we have $h_k (t) \ge \eta (t) / \mu$ for all $t \in (t_s, t_{s+1})$.
Thus, in view of~\eqref{pl3.5.5}, the estimate
$$
	\int\limits_{t_s}^{t_k}
	(h_k (t) t)^{- 1 / 4}
	dt
	\le
	\frac{
		\mu^{1 / 4}
	}{
		1 - \mu^{- 1 / 16}
	}
	\int\limits_{t_s}^{t_{s+1}}
	(\eta (t) t)^{- 1 / 4}
	dt
$$
is valid.
Since $h_s = h_k$ on the interval $[1, t_s]$, adding this with~\eqref{pl3.5.3}, we again arrive at~\eqref{pl3.5.1}.

By~\eqref{l3.5.2},
$$
	\lim_{k \to \infty}
	\gamma_k
	=
	\infty;
$$
therefore, the sequence $h_k$, $k = 1,2,\ldots$, stabilizes at $k \to \infty$ on any compact set belonging to the set $[1, \infty)$.
In particular, for all $t \in [1, \infty)$ there exists a limit
$$
	\tilde h (t)
	=
	\lim_{k \to \infty}
	h_k (t).
$$
Let us extend the function $\tilde h$ to the entire real axis, putting
$$
	\tilde h (t)
	=
	\min 
	\left\{
		\tilde h (1),
		\inf_{
			(t, 1]
		}
		H
	\right\}
$$
for $t \in (0, 1)$ and $\tilde h (t) = 0$ for $t \in (-\infty, 0]$.
To complete the proof of Lemma~5, it remasins to take
$$
	h (t)
	=
	\int\limits_{-\infty}^\infty
	\omega_{1/2} (t - \tau)
	\tilde h 
	\left(
		\tau - \frac{1}{2} 
	\right)
	d\tau,
	\quad
	t > 0.
$$

We prove Theorem~1 in several steps.
Let us agree to denote by $C$ positive, possibly different, constants depending only on $\sigma$ and $\theta$.

{\bf Step 1.}
Applying Lemma~5 with $\mu = \theta^{1/2}$ and
$$
	H (t)
	=
	\sup_{
		(t / \mu, t \mu)
	}
	g,
$$
we have a non-decreasing infinitely smooth function $h : (0, \infty) \to (0, \infty)$ such that
\begin{equation}
	h (4 t) \le C h (t)
	\label{3.14}
\end{equation}
for all $t \ge 2$, condition~\eqref{l3.5.3} is valid and, moreover,
\begin{equation}
	\inf_{
		x \in \Omega_{r, \sigma},
		\,
		\zeta \in (\theta^{-1/2} t, \theta^{1/2} t)
	}
	f (x, |x| \zeta)
	\ge
	q (r)
	h (t)
	\label{3.1}
\end{equation}
for all $t > 0$ and $r > r_*$.

{\bf Step 2.}
We take a continuous function
$p : [0, \infty) \to [0, \infty)$ and a real number $R_* > r_*$, 
satisfying the following conditions:
$$
	p (r) = 0 
	\quad
	\mbox{for all } 0 \le r \le R_*;
$$
\begin{equation}
	p (r)
	=
	e^{
		- 1 / (r - R_*) + 1
	}
	p (R_* + 1)
	>
	0
	\quad
	\mbox{for all } R_* < r \le R_* + 1;
	\label{3.4}
\end{equation}
\begin{equation}
	\inf_{
		x \in \Omega_{r, \sigma^{1/4}},
		\,
		\zeta \in (\theta^{-1/2} t, \theta^{1/2} t)
	}
	f (x, |x| \zeta)
	\ge
	p (r)
	h (t)
	\quad
	\mbox{for all } t > 0, r > R_*;
	\label{3.5}
\end{equation}
\begin{equation}
	r p (r) \le 1
	\quad
	\mbox{for all } r \ge 0;
	\label{3.7}
\end{equation}
\begin{equation}
	\int\limits_{r_*}^\infty
	p (r)
	dr
	=
	\infty.
	\label{3.6}
\end{equation}

Let us show that such a function $p$ and a real number $R_*$ exist, and $R_*$ can be made arbitrarily large. We put
$$
	Q_\delta (r)
	=
	\int\limits_{-\infty}^\infty
	\omega_\delta (r - \rho)
	Q (\rho)
	d\rho,
	\quad
	\delta > 0,
$$
where
$$
	Q (\rho)
	=
	\left\{
		\begin{aligned}
			&
			0,
			&
			&
			\rho \in (-\infty, r_*],
			\\
			&
			\sup_{
				(\sigma^{-1/4} \rho, \sigma^{1/4} \rho)
				\cap
				(r_*, \infty)
			}
			q,
			&
			&
			\rho \in (r_*, \infty).
		\end{aligned}
	\right.
$$
Taking into account~\eqref{3.1}, we have
\begin{equation}
	\inf_{
		x \in \Omega_{r, \sigma^{1/2}},
		\,
		\zeta \in (\theta^{-1/2} t, \theta^{1/2} t)
	}
	f (x, |x| \zeta)
	\ge
	Q_\delta (r)
	h (t)
	\label{3.2}
\end{equation}
for all $t > 0$ and $r > r_*$, where $\delta > 0$ is a sufficiently small real number.
At the same time, for any $\delta > 0$ there is a real number $r_\delta > r_*$ such that
$$
	\inf_{
		(\sigma^{-1/8} r, \sigma^{1/8}r)
	}
	Q_\delta
	\ge
	q (r)
$$
for all $r \in [r_\delta, \infty)$.
Denote
$$
	M_\delta (r)
	=
	\min 
	\left\{
		Q_\delta (r),
		\frac{1}{r}
	\right\}.
$$
We show that
\begin{equation}
	\int\limits_{r_*}^\infty
	M_\delta (r)
	dr
	=
	\infty
	\label{3.3}
\end{equation}
for any $\delta > 0$.
Indeed, let there be a real number $r > 0$ in every neighborhood of infinity such that
$q (r) \ge 1 / r$.
Then
$M_\delta (\rho) \ge \sigma^{-1/8} / \rho$
for all
$\rho \in (\sigma^{-1/8} r, \sigma^{1/8} r)$;
therefore,
$$
	\int\limits_{
		\sigma^{-1/8} r
	}^{
		\sigma^{1/8} r
	}
	M_\delta (\rho)
	d\rho
	=
	\frac{1}{4}
	\sigma^{-1/8}
	\ln \sigma.
$$
According to the necessary criterion for the convergence of an improper integral, this proves~\eqref{3.3}.
Suppose now that $q (r) < 1 / r$ for all $r > 0$ in some neighborhood of infinity.
In this case, we obtain
$$
	M_\delta (r)
	\ge
	\min 
	\left\{
		q (r),
		\frac{1}{r}
	\right\}
	=
	q (r)
$$
for all sufficiently large $r > 0$. Hence,~\eqref{3.3} follows from the second condition of~\eqref{t2.1.1}.

Thus, it remains to take
$$
	p (r)
	=
	\left\{
		\begin{aligned}
			&
			0,
			&
			&
			r \in [0, R_*],
			\\
			&
			e^{
				- 1 / (r - R_*) + 1
			}
			M_\delta (R_* + 1),
			&
			&
			r \in (R_*, R_* + 1],
			\\
			&
			M_\delta (r),
			&
			&
			r \in (R_* + 1, \infty),
		\end{aligned}
	\right.
$$
where the real numbers $R_* > r_*$ and $\delta > 0$ are chosen so that
${M_\delta (R_* + 1)} > 0$ and, moreover,~\eqref{3.5} holds.
The existence of such real numbers follows from~\eqref{3.2} and~\eqref{3.3}.
It can be seen that, in the case where $R_*$ is large enough, the relation~\eqref{3.2} implies~\eqref{3.5}.
Indeed, if $R_* > 1 / (\sigma^{1/4} - 1)$, then for all $r \in (R_*, R_* + 1)$ the inclusion
$
	\Omega_{r, \sigma^{1/4}} 
	\subset 
	\Omega_{R_* + 1, \sigma^{1/2}}
$
is valid; therefore,
$$
	\inf_{
		x \in \Omega_{r, \sigma^{1/4}},
		\,
		\zeta \in (\theta^{-1/2} t, \theta^{1/2} t)
	}
	f (x, |x| \zeta)
	\ge
	\inf_{
		x \in \Omega_{R_* + 1, \sigma^{1/2}},
		\,
		\zeta \in (\theta^{-1/2} t, \theta^{1/2} t)
	}
	f (x, |x| \zeta)
$$
for all $r \in (R_*, R_* + 1)$.
Since the function $p$ does not exceed $Q_\delta (R_* + 1)$ on the interval $(R_*, R_* + 1)$, this allows us to assert that~\eqref{3.5} holds for all $ r \in (R_*, R_* + 1)$.
Finally, if $r \in [R_* + 1, \infty)$, then the inequality~\eqref{3.2} implies~\eqref{3.5}, because for these $r$ the value of the function $p$ does not exceed $Q_\delta$.
Note also that, in view of~\eqref{3.3}, in any neighborhood of infinity there is a point at which the function $M_\delta$ is positive.

{\bf Step 3.}
Consider the Cauchy problem for the ordinary differential equation
\begin{equation}
	\frac{
		d^2 w
	}{
		dr^2
	}
	\left(
		\frac{
			d w
		}{
			dr
		}
	\right)^2
	=
	\frac{1}{2}
	p (r)
	h 
	\left(
		\frac{w}{r}
	\right)
	\quad
	\mbox{for } r > 0,
	\quad
	w (0) = \varepsilon,
	\quad
	\frac{
		d w
	}{
		d r
	}
	(0)
	=
	0,
	\label{3.8}
\end{equation}
where $\varepsilon > 0$ is some real number.
It is obvious that a local solution of this problem exists. 
In particular, on the interval $[0, R_*]$, this solution is the function $w = \varepsilon$.
Let 
$$
	R_{max}
	=
	\sup {\mathcal R},
$$
where ${\mathcal R}$ is the set of real numbers $R > 0$ such that the solution of problem~\eqref{3.8} exists on the interval $[0, R)$.
Reducing~\eqref{3.8} to the integral equation
\begin{equation}
	w (r)
	=
	\varepsilon
	+
	\int\limits_0^r
	\left(
		\frac{3}{2}
		\int\limits_0^\rho
		p (\xi) 
		h 
		\left(
			\frac{w}{\xi}
		\right)
		d\xi
	\right)^{1/3}
	d\rho,
	\label{3.9}
\end{equation}
by standard reasoning, for example, by the contraction mapping method, one can show that $R_{max} > R_*$ and on the interval $[0, R_{max})$ the solution of problem~\eqref{3.8} is unique.
Also it is easy to show that $w \in C^2 ([0, R_{max}))$.
Indeed, $w$ is an infinitely smooth function on the set $[0, R_*)$,
because on this set $w$ is a constant.
At the same time, differentiating~\eqref{3.9}, we have
$$
	\frac{d w (r)}{dr}
	=
	\left(
		\frac{3}{2}
		\int\limits_0^r
		p (\xi) 
		h 
		\left(
			\frac{w}{\xi}
		\right)
		d\xi
	\right)^{1/3},
	\quad
	\frac{d^2 w (r)}{dr^2}
	=
	\frac{1}{2}
	p (r) 
	h 
	\left(
		\frac{w}{r}
	\right)
	\left(
		\frac{3}{2}
		\int\limits_0^r
		p (\xi) 
		h 
		\left(
			\frac{w}{\xi}
		\right)
		d\xi
	\right)^{-1/3}
$$
for all $r \in (R_*, R_{max})$.
Since
$$
	\int\limits_0^r
	p (\xi) 
	h 
	\left(
		\frac{w}{\xi}
	\right)
	d\xi
	>
	0
$$
for all $r \in (R_*, R_{max})$,
the function $w$ is twice continuously differentiable on the set $(R_*, R_{max})$.
According to~\eqref{3.4}, the relation
$$
	\int\limits_0^r
	p (\xi) 
	h 
	\left(
		\frac{w}{\xi}
	\right)
	d\xi
	=
	(r - R_*)^2
	e^{- 1 / (r - R_*) + 1}
	p (R_* + 1)
	h 
	\left(
		\frac{
			\varepsilon
		}{
			R_*
		}
	\right)
	+
	\bar{\bar{o}} 
	\left(
		(r - R_*)^2
		e^{- 1 / (r - R_*) + 1}
	\right)
$$
is valid as $r \to R_* + 0$; 
therefore the function $w$ is twice differentiable at the point $R_*$ with
$$
	\frac{d w (R_*)}{dr}
	=
	\frac{d^2 w (R_*)}{dr^2}
	=
	0.
$$
Thus, it remains to note that
$$
	\lim_{r \to R_*}
	\frac{d w (r)}{dr}
	=
	\lim_{r \to R_*}
	\frac{d^2 w (r)}{dr^2}
	=
	0.
$$

Now we establish the validity of the inequality
\begin{equation}
	R_{max} < \infty.
	\label{3.10}
\end{equation}
Assume the contrary, let the solution of the problem~\eqref{3.8} exist on the entire interval $[0, \infty)$.
Putting
$$
	\frac{w (r) - \varepsilon}{r}
	=
	v (r)
$$
in equation~\eqref{3.9}, we obtain
\begin{equation}
	v (r)
	=
	\frac{1}{r}
	\int\limits_0^r
	\left(
		\frac{3}{2}
		\int\limits_0^\rho
		p (\xi) 
		h 
		\left(
			v
			+
			\frac{\varepsilon}{\xi}
		\right)
		d\xi
	\right)^{1/3}
	d\rho.
	\label{3.11}
\end{equation}
Differentiating~\eqref{3.11}, we have
\begin{equation}
	\frac{d v (r)}{dr}
	=
	\frac{1}{r^2}
	\left(
		r \varphi (r)
		-
		\int\limits_0^r
		\varphi (\rho)
		d\rho
	\right)
	\label{3.12}
\end{equation}
for all $r > 0$, where
\begin{equation}
	\varphi (\rho)
	=
	\left(
		\frac{3}{2}
		\int\limits_0^\rho
		p (\xi) 
		h 
		\left(
			v
			+
			\frac{\varepsilon}{\xi}
		\right)
		d\xi
	\right)^{1/3}.
	\label{3.17}
\end{equation}
Note that on the set $[0, R_*]$ the functions $v$ and $\varphi$ are identically equal to zero.
Since $\varphi$ is non-decreasing on $(0, \infty)$, the right-hand side of the formula~\eqref{3.12} is non-negative; therefore, $v$ is a non-decreasing function on $[0, \infty)$.
It can also be shown that
\begin{equation}
	\lim_{r \to \infty}
	v (r)
	=
	\infty.
	\label{3.15}
\end{equation}
Indeed, from~\eqref{3.6} and the fact that the function $h$ is positive, it follows that $v (\rho_*) > 0$ for some real number $\rho_* > R_*$.
Thus,
$$
	v (r)
	\ge
	\frac{1}{r}
	\int\limits_{r/2}^r
	\left(
		\frac{3}{2}
		\int\limits_{\rho_*}^\rho
		p (\xi) 
		h 
		\left(
			v
			+
			\frac{\varepsilon}{\xi}
		\right)
		d\xi
	\right)^{1/3}
	d\rho
	\ge
	\frac{1}{2}
	h^{1/3} (v (\rho_*))
	\left(
		\frac{3}{2}
		\int\limits_{\rho_*}^{r/2}
		p (\xi) 
		d\xi
	\right)^{1/3}
	\to
	\infty
$$
as $r \to \infty$.
Take a real number $r_0 > r_*$ such that $\varepsilon / r_0 \le 1$, $v (r_0) \ge 4$ and
\begin{equation}
	\frac{
		r p(r)
	}{
		\int\limits_{r_*}^r
		p (\xi) 
		d\xi
	}
	\le
	1
	\label{3.13}
\end{equation}
for all $r \ge r_0$.
In view of~\eqref{3.15}, \eqref{3.7}, and~\eqref{3.6} such a real number $r_0$ exists.
We put
$$
	r_i 
	=
	\sup
	\{
		r \in (r_{i-1}, 2 r_{i-1})
		:
		v (r) \le 2 v (r_{i-1})
	\},
	\quad
	i = 1,2,\ldots.
$$
In is easy to see that $r_i \to \infty$ as $i \to \infty$; otherwise~\eqref{3.10} holds.

\textbf{Lemma~6.} {\it 
Let $v (r_{i+1}) \ge 2 v (r_{i-1})$ for some $i \ge 1$.
Then
\begin{equation}
	\int\limits_{
		v (r_{i-1})
	}^{
		v (r_{i+1})
	}
	(h (t) t)^{- 1 / 4}
	dt
	\ge
	\frac{
		C
		\int\limits_{
			r_{i-1}
		}^{
			r_{i+1}
		}
		p (\xi)
		d\xi
	}{
		\left(
			\int\limits_{
				r_*
			}^{
				r_{i+1}
			}
			p (\xi)
			d\xi
		\right)^{3/4}
	}.
	\label{l3.6.1}
\end{equation}
}

{\bf Proof.} 
By~\eqref{3.11}, the inequality
$$
	v (r_{i+1})
	\ge
	\frac{
		h^{1/3} (v (r_{i-1}))
	}{
		r_{i+1}
	}
	\int\limits_{
		r_{i-1}
	}^{
		r_{i+1}
	}
	\left(
		\frac{3}{2}
		\int\limits_{
			r_{i-1}
		}^\rho
		p (\xi) 
		d\xi
	\right)^{1/3}
	d\rho
$$
iv valid. Combining this with the estimate
$$
	\int\limits_{
		r_{i-1}
	}^{
		r_{i+1}
	}
	\left(
		\int\limits_{
			r_{i-1}
		}^\rho
		p (\xi) 
		d\xi
	\right)^{1/3}
	d\rho
	\ge
	C
	\int\limits_{
		r_{i-1}
	}^{
		r_{i+1}
	}
	((r_{i+1}	- \xi) p (\xi))^{1/3}
	d\xi
$$
which follows from Lemma~2, we obtain
\begin{equation}
	\frac{
		v (r_{i+1})
	}{
		h^{1/3} (v (r_{i-1}))
	}
	\ge
	\frac{
		C
	}{
		r_{i+1}
	}
	\int\limits_{
		r_{i-1}
	}^{
		r_{i+1}
	}
	((r_{i+1}	- \xi) p (\xi))^{1/3}
	d\xi.
	\label{pl3.6.1}
\end{equation}
Applying further Lemma~3, we arrive at the relation
\begin{equation}
	\left(
		\int\limits_{
			r_{i-1}
		}^{
			r_{i+1}
		}
		((r_{i+1}	- \xi) p (\xi))^{1/3}
		d\xi
	\right)^{3/4}
	\ge
	C
	\int\limits_{
		r_{i-1}
	}^{
		r_{i+1}
	}
	((r_{i+1}	- \xi) p (\xi))^{1/3}
	\varkappa^{-1/4} (\xi)
	d\xi,
	\label{pl3.6.2}
\end{equation}
where
$$
	\varkappa (\xi)
	=
	\int\limits_\xi^{r_{i+1}}
	((r_{i+1}	- \zeta) p (\zeta))^{1/3}
	d\zeta.
$$
It is easy to see that
\begin{align*}
	(r_{i+1}	- \xi)^{1/3}
	\varkappa^{-1/4} (\xi)
	&
	=
	\left(
		\frac{
			1
		}{
			(r_{i+1}	- \xi)^{4/3}
		}
		\int\limits_\xi^{r_{i+1}}
		((r_{i+1}	- \zeta) p (\zeta))^{1/3}
		d\zeta
	\right)^{-1/4}
	\\
	&
	\ge
	\left(
		\frac{
			1
		}{
			r_{i+1}	- \xi
		}
		\int\limits_\xi^{r_{i+1}}
		p^{1/3} (\zeta)
		d\zeta
	\right)^{-1/4}
\end{align*}
for all $\xi \in (r_{i-1}, r_{i+1})$.
At the same time, in view of~\eqref{3.13}, we have
$$
	p^{1/3} (\zeta)
	\le
	\zeta^{-1/3}
	\left(
		\int\limits_{
			r_*
		}^{
			\zeta
		}
		p (\rho)
		d\rho
	\right)^{1/3}
	\le
	\xi^{-1/3}
	\left(
		\int\limits_{
			r_*
		}^{
			r_{i+1}
		}
		p (\rho)
		d\rho
	\right)^{1/3}
$$
for all $\zeta \in (\xi, r_{i+1})$. Thus,
$$
	(r_{i+1}	- \xi)^{1/3}
	\varkappa^{-1/4} (\xi)
	\ge
	\xi^{1/12}
	\left(
		\int\limits_{
			r_*
		}^{
			r_{i+1}
		}
		p (\rho)
		d\rho
	\right)^{-1/12}
$$
for all $\xi \in (r_{i-1}, r_{i+1})$ and the inequality~\eqref{pl3.6.2} allows us to assert that
$$
	\left(
		\int\limits_{
			r_{i-1}
		}^{
			r_{i+1}
		}
		((r_{i+1}	- \xi) p (\xi))^{1/3}
		d\xi
	\right)^{3/4}
	\ge
	\frac{
		C
		\int\limits_{
			r_{i-1}
		}^{
			r_{i+1}
		}
		\xi^{1/12}
		p^{1/3} (\xi)
		d\xi
	}{
		\left(
			\int\limits_{
				r_*
			}^{
				r_{i+1}
			}
			p (\rho)
			d\rho
		\right)^{1/12}
	}.
$$
Combining the last expression with formula~\eqref{pl3.6.1}, we obtain
\begin{equation}
	\frac{
		v^{3/4} (r_{i+1})
	}{
		h^{1/4} (v (r_{i-1}))
	}
	\ge
	\frac{
		C
		\int\limits_{
			r_{i-1}
		}^{
			r_{i+1}
		}
		\xi^{-2/3}
		p^{1/3} (\xi)
		d\xi
	}{
		\left(
			\int\limits_{
				r_*
			}^{
				r_{i+1}
			}
			p (\rho)
			d\rho
		\right)^{1/12}
	}.
	\label{pl3.6.3}
\end{equation}
According to~\eqref{3.13}, 
$$
	\xi^{-2/3}
	p^{1/3} (\xi)
	=
	(p (\xi) \xi)^{-2/3}
	p (\xi)
	\ge
	\frac{
		p (\xi)
	}{
		\left(
			\int\limits_{
				r_*
			}^{
				\xi
			}
			p (\rho)
			d\rho
		\right)^{2/3}
	}
	\ge
	\frac{
		p (\xi)
	}{
		\left(
			\int\limits_{
				r_*
			}^{
				r_{i+1}
			}
			p (\rho)
			d\rho
		\right)^{2/3}
	},
$$
for all $\xi \in (r_{i-1}, r_{i+1})$; therefore,~\eqref{pl3.6.3} leads to the estimate
\begin{equation}
	\frac{
		v^{3/4} (r_{i+1})
	}{
		h^{1/4} (v (r_{i-1}))
	}
	\ge
	\frac{
		C
		\int\limits_{
			r_{i-1}
		}^{
			r_{i+1}
		}
		p (\xi)
		d\xi
	}{
		\left(
			\int\limits_{
				r_*
			}^{
				r_{i+1}
			}
			p (\rho)
			d\rho
		\right)^{3/4}
	}.
	\label{pl3.6.4}
\end{equation}
Considering~\eqref{3.14} and the fact that 
$
	v (r_{i+1}) - v (r_{i-1}) \ge v (r_{i+1}) / 2,	
$ 
we have
$$
	\int\limits_{
		v (r_{i-1})
	}^{
		v (r_{i+1})
	}
	(h (t) t)^{-1/4}
	dt
	\ge
	\frac{
		v (r_{i+1}) - v (r_{i-1})
	}{
		(h (v (r_{i+1})) v (r_{i+1}))^{1/4}
	}
	\ge
	\frac{
		v^{3/4} (r_{i+1})
	}{
		2 h^{1/4} (v (r_{i+1}))
	}
	\ge
	\frac{
		C
		v^{3/4} (r_{i+1})
	}{
		h^{1/4} (v (r_{i-1}))
	}.
$$
Thus,~\eqref{pl3.6.4} implies~\eqref{l3.6.1}.
Lemma~6 is proved.

{\bf Lemma~7.} {\it
Let $v (r_{i+1}) < 2 v (r_{i-1})$ for some $i \ge 1$.
Then
\begin{equation}
	\int\limits_{
		v (r_{i-1})
	}^{
		v (r_{i+1})
	}
	\frac{
		dt
	}{
		h^{1/3} (t)
	}
	\ge
	\frac{
		C
		\int\limits_{
			r_{i-1}
		}^{
			r_i
		}
		p (\xi)
		d\xi
	}{
		\left(
			\int\limits_{
				r_*
			}^{
				r_i
			}
			p (\rho)
			d\rho
		\right)^{2/3}
	}.
	\label{l3.7.1}
\end{equation}
}

{\bf Proof.} 
It follows from the conditions of the lemma that $r_i = 2 r_{i-1}$ and $r_{i+1} = 2 r_i$.
By~\eqref{3.12}, we get
\begin{equation}
	v (r_{i+1}) - v (r_{i-1})
	=
	\int\limits_{
		r_{i-1}
	}^{
		r_{i+1}
	}
	v' (r)
	dr
	\ge
	\frac{1}{
		r_{i+1}^2
	}
	\int\limits_{
		r_{i-1}
	}^{
		r_{i+1}
	}
	\psi (r)
	dr,
	\label{pl3.7.1}
\end{equation}
where
$$
	\psi (r)
	=
	r \varphi (r)
	-
	\int\limits_0^r
	\varphi (\rho)
	d\rho.
$$
By direct differentiation, one can verify that
$
	\psi' (r) 
	= 
	r \varphi' (r)
	\ge
	0
$
for all $r > 0$.
According to the Newton-Leibniz formula
$$
	\psi (r)
	=
	\int\limits_0^r
	\xi
	\varphi' (\xi)
	d\xi
	\ge
	\int\limits_{r_{i-1}}^r
	\xi 
	\varphi' (\xi)
	d\xi
$$
for all $r \in (r_{i-1}, r_{i+1})$;
therefore,~\eqref{pl3.7.1} implies the estimate
$$
	v (r_{i+1}) - v (r_{i-1})
	\ge
	\frac{
		1
	}{
		r_{i+1}^2
	}
	\int\limits_{
		r_{i-1}
	}^{
		r_{i+1}
	}
	\int\limits_{r_{i-1}}^r
	\xi 
	\varphi' (\xi)
	d\xi
	dr
	\ge
	\frac{
		C
	}{
		r_{i+1}
	}
	\int\limits_{
		r_{i-1}
	}^{
		r_{i+1}
	}
	\int\limits_{r_{i-1}}^r
	\varphi' (\xi)
	d\xi
	dr,
$$
combining which with the elementary inequality
$$
	\int\limits_{
		r_{i-1}
	}^{
		r_{i+1}
	}
	\int\limits_{r_{i-1}}^r
	\varphi' (\xi)
	d\xi
	dr
	\ge
	\int\limits_{
		r_i
	}^{
		r_{i+1}
	}
	\int\limits_{r_{i-1}}^r
	\varphi' (\xi)
	d\xi
	dr
	\ge
	(r_{i+1} - r_i)
	\int\limits_{r_{i-1}}^{r_i}
	\varphi' (\xi)
	d\xi,
$$
we have
\begin{equation}
	v (r_{i+1}) - v (r_{i-1})
	\ge
	C
	\int\limits_{r_{i-1}}^{r_i}
	\varphi' (\xi)
	d\xi.
	\label{pl3.7.2}
\end{equation}
At the same time, differentiating~\eqref{3.17}, we obtain
$$
	\varphi' (\xi)
	=
	\frac{
		2^{-1/3}
		3^{-2/3}
		p (\xi) 
		h 
		\left(
			v (\xi)
			+
			\varepsilon / \xi
		\right)
	}{
		\left(
			\int\limits_{r_*}^\xi
			p (\rho) 
			h 
			\left(
				v
				+
				\varepsilon / \rho
			\right)
			d\rho
		\right)^{2/3}
	}
	\ge
	\frac{
		2^{-1/3}
		3^{-2/3}
		p (\xi) 
		h (v (r_{i-1}))
	}{
		\left(
			\int\limits_{
				r_*
			}^{
				r_i
			}
			p (\rho) 
			d\rho
		\right)^{2/3}
		h^{2/3} (2 v (r_i))
	}
$$
for all $\xi \in (r_{i-1}, r_i)$. 
Since $h (2 v (r_i)) \le C h (v (r_{i-1}))$ in view of the condition~\eqref{3.14}, this allows us to assert that
$$
	\varphi' (\xi)
	\ge
	\frac{
		C
		p (\xi) 
		h^{1/3} (v (r_{i-1}))
	}{
		\left(
			\int\limits_{
				r_*
			}^{
				r_i
			}
			p (\rho) 
			d\rho
		\right)^{2/3}
	}
$$
for all $\xi \in (r_{i-1}, r_i)$. Thus, the formula~\eqref{pl3.7.2} implies the estimate
$$
	\frac{
		v (r_{i+1}) - v (r_{i-1})
	}{
		h^{1/3} (v (r_{i-1}))
	}
	\ge
	\frac{
		C
		\int\limits_{r_{i-1}}^{r_i}
		p (\xi)
		d\xi
	}{
		\left(
			\int\limits_{
				r_*
			}^{
				r_i
			}
			p (\rho) 
			d\rho
		\right)^{2/3}
	},
$$
combining which with the inequality
$$
	\int\limits_{
		v (r_{i-1})
	}^{
		v (r_{i+1})
	}
	\frac{
		dt
	}{
		h^{1/3} (t)
	}
	\ge
	\frac{
		C (v (r_{i+1}) - v (r_{i-1}))
	}{
		h^{1/3} (v (r_{i-1}))
	},
$$
we arrive at~\eqref{l3.7.1}. Lemma~7 is proved.

For a natural number $m$, we denote by $\Xi_{1, m}$ the set of integers ${1 \le i \le m}$ satisfying the condition $v (r_{i+1}) \ge 2 v (r_ {i-1})$. Also let
$
	\Xi_{2, m} 
	= \{ 1, \ldots, m \} 
	\setminus 
	\Xi_{1, m}.
$

In view of~\eqref{3.15} and~\eqref{3.6}, there is a natural number $l$ such that
$v (r_m) \ge 2 v (r_0)$ 
and
\begin{equation}
	\int\limits_{
		r_0
	}^{
		r_m
	}
	p (\xi)
	d\xi
	\ge
	\frac{1}{2}
	\int\limits_{
		r_*
	}^{
		r_m
	}
	p (\xi)
	d\xi
	\label{3.20}
\end{equation}
for all $m \ge l$.
Let us show that
\begin{equation}
	\int\limits_{
		v (r_0)
	}^{
		v (r_{m + 1})
	}
	(h (t) t)^{- 1 / 4}
	dt
	\ge
	C
	\left(
		\int\limits_{
			r_*
		}^{
			r_{m + 1}
		}
		p (\xi)
		d\xi
	\right)^{1/4}
	\label{3.16}
\end{equation}
for all $m \ge l$.
At first, we assume that
\begin{equation}
	\sum_{
		i \in \Xi_{1, m}
	}
	\,
	\int\limits_{
		r_{i-1}
	}^{
		r_i
	}
	p (\xi)
	d\xi
	\ge
	\frac{1}{2}
	\int\limits_{r_0}^{r_m}
	p (\xi)
	d\xi.
	\label{3.18}
\end{equation}
Summing estimate~\eqref{l3.6.1} of Lemma~6 over all $i \in \Xi_{1, m}$, we have
$$
	\int\limits_{
		v (r_0)
	}^{
		v (r_{m + 1})
	}
	(h (t) t)^{- 1 / 4}
	dt
	\ge
	\frac{
		C
		\int\limits_{
			r_0
		}^{
			r_m
		}
		p (\xi)
		d\xi
	}{
		\left(
			\int\limits_{
				r_*
			}^{
				r_{m + 1}
			}
			p (\xi)
			d\xi
		\right)^{3/4}
	},
$$
whence in view of~\eqref{3.20} and the inequality
\begin{equation}
	\int\limits_{
		r_*
	}^{
		r_m
	}
	p (\xi)
	d\xi
	\ge
	C
	\int\limits_{
		r_*
	}^{
		r_{m + 1}
	}
	p (\xi)
	d\xi,
	\label{3.19}
\end{equation}
which follows from Lemma~4, we immediately arrive at~\eqref{3.16}.
Now assume that~\eqref{3.18} is not valid. In this case, it obvious that
$$
	\sum_{
		i \in \Xi_{2, m}
	}
	\,
	\int\limits_{
		r_{i-1}
	}^{
		r_i
	}
	p (\xi)
	d\xi
	\ge
	\frac{1}{2}
	\int\limits_{r_0}^{r_m}
	p (\xi)
	d\xi.
$$
Thus, summing estimate~\eqref{l3.7.1} of Lemma~7 over all $i \in \Xi_{2, m}$, we arrive at the relation
$$
	\int\limits_{
		v (r_0)
	}^{
		v (r_{m + 1})
	}
	\frac{
		dt
	}{
		h^{1/3} (t)
	}
	\ge
	\frac{
		C
		\int\limits_{
			r_0
		}^{
			r_m
		}
		p (\xi)
		d\xi
	}{
		\left(
			\int\limits_{
				r_*
			}^{
				r_m
			}
			p (\rho)
			d\rho
		\right)^{2/3}
	},
$$
combining which with~\eqref{3.20} and~\eqref{3.19}, we obtain
\begin{equation}
	\int\limits_{
		v (r_0)
	}^{
		v (r_{m + 1})
	}
	\frac{
		dt
	}{
		h^{1/3} (t)
	}
	\ge
	C
	\left(
		\int\limits_{
			r_*
		}^{
			r_{m + 1}
		}
		p (\rho)
		d\rho
	\right)^{1/3}.
	\label{3.21}
\end{equation}
Since $h$ is a non-decreasing function satisfying condition~\eqref{3.14}, Lemma~1 yields
$$
	\left(
		\int\limits_{
			v (r_0)
		}^{
			v (r_{m + 1})
		}
		(h (t) t)^{- 1 / 4}
		dt
	\right)^{4/3}
	\ge
	C
	\int\limits_{
		v (r_0)
	}^{
		v (r_{m + 1})
	}
	\frac{
		dt
	}{
		h^{1/3} (t)
	}.
$$
Combining the last estimate with~\eqref{3.21}, we again obtain~\eqref{3.16}.

Passing in~\eqref{3.16} to the limit as $m \to \infty$, we obtain a contradiction with the conditions~\eqref{l3.5.3} and~\eqref{3.6}. This proves the validity of~\eqref{3.10}.
Note that, in view of~\eqref{3.10}, the relation
\begin{equation}
	\lim_{r \to R_{\max} - 0} w (r) = \infty;
	\label{3.22}
\end{equation}
otherwise the solution of the problem~\eqref{3.8} exists on the interval $(0, R)$ for some $R > R_{max}$.

{\bf Step 4.}
Below we will assume that $u$ is a non-negative solution~\eqref{1.1}.
We need the comparison principle in the following simple form.

\textbf{Lemma~8.} {\it 
Let
$U \in C^2 (B_R \cap \Omega) \cap C (\overline{B_R \cap \Omega})$ 
be a classical solution of the inequality
$$
	\Delta_\infty U \le {\mathcal F} (x, U)
	\quad
	\mbox{in } B_R \cap \Omega
$$
non-negative on the set $B_R \cap \Omega$ and such that
\begin{equation}
	\left.
		U
	\right|_{
		\partial (B_R \cap \Omega)
	}
	\ge
	\left.
		u
	\right|_{
		\partial (B_R \cap \Omega)
	},
	\label{l3.8.2}
\end{equation}
where $R > 0$ is a real number and ${\mathcal F}$ is some function non-decreasing with respect to the last argument with
\begin{equation}
	{\mathcal F} (x, t) < f (x, t)
	\label{l3.8.3}
\end{equation}
for all $x \in B_R \cap \Omega$ and $t > 0$. Then
\begin{equation}
	U \ge u
	\quad
	\mbox{in } \overline{B_R \cap \Omega}.
	\label{l3.8.4}
\end{equation}
}

\textbf{Proof.}
Assume the contrary. Let
$$
	\sup_{
		\overline{B_R \cap \Omega}
	} 
	(u - U)
	>
	0.
$$
The function $u - U$ is upper semicontinuous on the set $\overline{B_R \cap \Omega}$; therefore, there is a point $x \in \overline{B_R \cap \Omega}$ such that
$$
	u (x) - U (x)
	=
	\sup_{
		\overline{B_R \cap \Omega}
	} 
	(u - U).
$$
In view of~\eqref{l3.8.2}, one can also assert that $x \in B_R \cap \Omega$. Hence,
$$
	u (y) - u (x)
	\le
	U (y) - U (x)
	=
	\langle D U (x), y - x \rangle
	+
	\frac{1}{2}
	\langle D^2 U (x) (y - x), y - x \rangle
	+
	\bar{\bar{o}} (|y - x|^2)
$$
for all $y$ from some neighborhood of the point $x$, where $D U (x)$ is the gradient and $D^2 U (x)$ is the matrix of the second derivatives of the function $U$ at the point $x$ , which immediately implies that
$$
	(D U (x), D^2 U (x)) \in J_\Omega^{2,+} u (x).
$$
Substituting
$v = D U (x)$ and $A = D^2 U (x)$
in~\eqref{d1.3.1}, we have
$$
	\sum_{i,j=1}^n
	\frac{
		\partial^2 U (x)
	}{
		\partial x_i \partial x_j
	}
	\frac{
		\partial U (x)
	}{
		\partial x_i
	}
	\frac{
		\partial U (x)
	}{
		\partial x_j
	}
	\ge
	f (x, u (x))
	>
	{\mathcal F} (x, u (x)).
$$
Thus, taking into account the inequalities
$$
	{\mathcal F} (x, u (x))
	\ge
	{\mathcal F} (x, U (x))
	\ge
	\sum_{i,j=1}^n
	\frac{
		\partial^2 U (x)
	}{
		\partial x_i \partial x_j
	}
	\frac{
		\partial U (x)
	}{
		\partial x_i
	}
	\frac{
		\partial U (x)
	}{
		\partial x_j
	},
$$
we conclude that
$$
	\sum_{i,j=1}^n
	\frac{
		\partial^2 U (x)
	}{
		\partial x_i \partial x_j
	}
	\frac{
		\partial U (x)
	}{
		\partial x_i
	}
	\frac{
		\partial U (x)
	}{
		\partial x_j
	}
	>
	\sum_{i,j=1}^n
	\frac{
		\partial^2 U (x)
	}{
		\partial x_i \partial x_j
	}
	\frac{
		\partial U (x)
	}{
		\partial x_i
	}
	\frac{
		\partial U (x)
	}{
		\partial x_j
	}.
$$
This contradiction proves the lemma.

\textbf{Proof of Theorem~1.}
According to~\eqref{3.8}, the function
$
	U (x) = w (|x|)
$ 
is a classical solution of the equation
$$
	\Delta_\infty U
	=
	{\mathcal F} (x, U)
	\quad
	\mbox{in } B_{R_{max}}
$$
positive in the ball $B_{R_{max}}$, where
$$
	{\mathcal F} (x, t)
	=
	\left\{
		\begin{aligned}
			&
			\frac{1}{2}
			p (|x|)
			h 
			\left( 
				\frac{t}{|x|}
			\right),
			&
			&
			R_* < |x| < R_{max},
			\,
			t > 0,
			\\
			&
			0,
			&
			&
			|x| \le R_*, 
			\,
			t > 0.
		\end{aligned}
	\right.
$$
Since a solution of~\eqref{1.1} is upper semicontinuous on $\overline{\Omega}$, the function $u$ is bounded on the closure of $B_{R_{max}} \cap \Omega$; therefore, in accordance with~\eqref{3.22} there is a real number $R \in (R_*, R_{max})$ for which~\eqref{l3.8.2} holds.
Taking into account~\eqref{1.2} and~\eqref{3.5}, one can also claim that~\eqref{l3.8.3} holds.
Thus, by Lemma~8, relation~\eqref{l3.8.4} is valid. In its turn, this implies that
$$
	\left.
		u
	\right|_{
		B_{R_*}
		\cap
		\Omega
	}
	\le
	\varepsilon.
$$
Since $\varepsilon > 0$ can be taken arbitrarily small and $R_* > r_*$ can be taken arbitrarily large, we obviously obtain $u (x) = 0$ for all $x \in \Omega$.

{\bf Proof of Corallary~1.} 
We put
$f (x, t) = c(|x|) t^\lambda$, 
$g (t) = t^\lambda$, 
$q (r) = 4^{-\lambda} r^\lambda c (r / 2)$,
and
$\sigma = \theta = 2$ in relation~\eqref{1.5} and use Theorem~1.

\pagebreak

\begin{center}{REFERENCES}\end{center}{\small

\lit{1}{Crandall~M.G., Ishii~H., Lions~P.-L.}{User's guide to viscosity solutions of second order partial differential equations //~Bull. Amer. Math. Soc. (N.S.). 1992. V.~27. P.~1--67.}
\lit{2}{Lu~G., Wang~P.}{Inhomogeneous infinity laplace equation //~Adv. Math. 2008. V.~217. P.~1838--1868.}
\lit{3}{Astashova~I.V.}{Uniqueness of solutions to second order Emden-Fowler type equations with general power-law nonlinearity //~J. Math. Sci. 2021. V.~225. No. 5. P.~543--550.}
\lit{4}{Astashova~I.V.}{Asymptotic behavior of singular solutions of Emden-Fowler type equations //~Diff. Eq. 2019. V.~55. No. 5. P.~581--590.}
\lit{5}{Astashova~I.V.}{On asymptotic behavior of blow-up solutions to higher-order differential equations with general nonlinearity //~Pinelas S., Caraballo T., Kloeden P., Graef J. (eds) Differential and Difference Equations with Applications. ICDDEA 2017. Springer Proceedings in Mathematics \& Statistics. Springer Cham. 2018. V. 230. P. 1--12.}
\lit{6}{Baras~P., Pierre~M.}{Singularit\'es \'eliminables pour des \'equations semilin\'eaires //~Ann. Inst. Fourier. 1984. V.~34. P.~185--205.}
\lit{7}{Galakhov~E.I.}{Solvability of an elliptic equation with a gradient nonlinearity //~Diff. Eq. 2005. V. 41. No. 5. P. 693--702.}
\lit{8}{Galakhov~E.I.}{Some nonexistence results for qu asilinear elliptic problems //~J. Math. Anal. Appl. 2000. V~252. No. 1. P.~256--277.}
\lit{9}{Fino A.Z., Galakhov E.I., Salieva O.A.}{Nonexistence of global weak solutions for evolution equations with fractional Laplacian //~Math. Notes 2020. V. 108. No. 5--6. P. 877--883.}
\lit{10}{Keller~J.B.}{On solutions of $\Delta u = f (u)$ //~Comm. Pure Appl. Math. 1957. V.~10. P.~503--510.}
\lit{11}{Kondratiev V.A., Landis E.M.}{Qualitative properties of the
solutions of a second-order nonlinear equations //~Sb. Math. 1989. V. 63. No. 2. P. 337--350.}
\lit{12}{Kon'kov~A.A.}{On global solutions of the radial $p$-Laplace equation //~Nonlinear Anal. 2009. V.~70. P.~3437--3451.}
\lit{13}{Kon'kov A.A.}{On solutions of nonautonomous ordinary differential equations //~Izv. Math. 2001. V. 65. No. 2. P. 285--327.}
\lit{14}{Kon'kov A.A.}{On properties of solutions of a class of nonlinear ordinary differential equations //~J. Math. Sci. 2007. V. 143. No. 4. P. 3303--3321.}
\lit{15}{Korpusov M.O., Matveeva A.K.}{On critical exponents for weak solutions of the Cauchy problem for a non-linear equation of composite type //~Izv. Math. 2021. V. 85. No. 4. P. 705--744.}
\lit{16}{Korpusov M.O., Panin A.A.}{On the nonextendable solution and blow-up of the solution of the one-dimensional equation of ion-sound waves in a plasma //~Math. Notes 2017. V. 102. No. 3--4. P. 350--360.}
\lit{17}{Korpusov M.O., Shafir R.S.}{Blow-up of weak solutions of the Cauchy problem for (3+1)-dimensional equation of plasma drift waves //~Comput. Math. Math. Phys. 2022. V. 62. No. 1. P. 117--149.}
\lit{18}{E.~Mitidieri~E.,Pohozaev~S.I.}{A priori estimates and blow-up of solutions to nonlinear partial differential equations and inequalities //~Proc. V.A.~Steklov Inst. Math. 2021. V. 234. P. 3--283.}
\lit{19}{Osserman~R.}{On the inequality $\Delta u \ge f (u)$ //~Pacific J. Math. 1957. V.~7. P.~1641--1647.}
\lit{20}{Mi~L.}{Blow-up rates of large solutions for infinity Laplace equations //~Appl. Math. Comp. 2017. V.~298. P.~36--44.}
\lit{21}{Mohammed~A., Mohammed~S.}{Boundary blow-up solutions to degenerate elliptic equations with non-monotone inhomogeneous terms //~Nonlinear Anal. 2012. V.~75. P.~3249--3261.}
\lit{22}{Wan~H.}{The exact asymptotic behavior of boundary blow-up solutions to infinity Laplacian equations //~Z. Angew. Math. Phys. 2016. V.~67. No.~97. P.~1--14.}
}

\end{document}